\documentclass[12pt]{amsart}
\usepackage{amscd, amssymb}

\theoremstyle{plain}

\newtheorem*{theorem*}{Theorem}

\theoremstyle{definition}

\newtheorem*{definition*}{Definition}

\newtheorem*{remark*}{Remark}

\theoremstyle{remark}

\def\be{\begin{equation}}
\def\ee{\end{equation}}
\def\this{i.\ e.\ } 
\def\h{\hbar}
\def\p{\partial}
\def\w{\wedge}

\def\id{\operatorname{Id}}

\def\dim{\operatorname{dim}}

\def\tr{\operatorname{tr}}
\def\str{\operatorname{str}}
\def\ev{\operatorname{ev}}
\def\ft{\operatorname{ft}}

\newcommand{\CC}{{\mathbb C}}

\newcommand{\ZZ}{{\mathbb Z}}
\newcommand{\QQ}{{\mathbb Q}}

\newcommand{\lan}{\langle}
\newcommand{\ran}{\rangle}
\newcommand{\Gb}{{\mathfrak b}} 

\newcommand{\Gg}{{\mathfrak g}}
\newcommand{\Gh}{{\mathfrak h}}

\newcommand{\Gn}{{\mathfrak n}}

\renewcommand{\a}{\alpha}
\renewcommand{\b}{\beta}

\renewcommand{\d}{\delta}
\newcommand{\e}{\varepsilon}

\renewcommand{\l}{\lambda}

\newcommand{\s}{\sigma}

\renewcommand{\v}{\nu}
\newcommand{\x}{\chi}
\newcommand{\z}{\zeta}

\newcommand{\D}{\Delta}
\renewcommand{\L}{\Lambda}

\renewcommand{\S}{\Sigma}

\newcommand{\calg}{\mathcal{G}}
\newcommand{\calh}{\mathcal{H}}

\newcommand{\calo}{\mathcal{O}}

\newcommand{\calt}{\mathcal{T}}

\newcommand{\calx}{QM}
\newcommand{\1}{{\bf 1}}
\newcommand{\E}{{\mathcal E}}

\title[Quantum $K$-theory on flag manifolds]
{Quantum $K$-theory on flag manifolds, \\ 
finite-difference Toda lattices \\ and quantum groups}

\author{Alexander GIVENTAL}
\address{UC Berkeley}
\email{givental@math.berkeley.edu}
\author{Yuan-Pin LEE}
\address{UCLA}
\email{yplee@math.ucla.edu}

\thanks{The research was partially supported by the NSF grants DMS-9704774 
(A. G.) and DMS-0072547 (Y. L.) }
\date{August 15, 2001}

\begin{document}

\maketitle

\setcounter{section}{-1}

\section{Introduction} \label{s:0}

Let $X=\{ 0\subset \CC^1\subset \ldots \subset \CC^r\subset \CC^{r+1}\}$
be the manifold of complete flags in $\CC ^{r+1}$. It admits the Pl\"ucker
embedding into the product of projective spaces 
\[
  X \hookrightarrow \Pi:= 
	 \prod_{i=1}^{r} \CC P^{n_i-1}, \quad n_i :=\binom{r+1}{i}.
\]

Let $(x:y)$ be homogeneous coordinates on $\CC P^1$.
A degree $d$ holomorphic map $\CC P^1\to \CC P^N$ is uniquely 
determined, up to a constant scalar factor,  by $N+1$ relatively prime
degree $d$ binary forms $(f_0(x:y):...:f_N(x:y))$. Omitting the condition
that the forms are relatively prime we compactify the space of degree $d$
holomorphic maps $\CC P^1\to \CC P^N$ to a complex projective space of
dimension $(N+1)(d+1)-1$. We denote this compactification of the space
of maps by $\CC P^N_d$. This construction defines the compactification
$\Pi_d = \Pi_{i=1}^r \CC P^{n_i-1}_{d_i}$ 
of the space of degree $d=(d_1,...,d_r)$ maps from $\CC P^1$ to $\Pi $. 

Composing degree $d$ holomorphic maps from $\CC P^1$ to the 
flag manifold $X$ with the Pl\"ucker embedding, we embed the space of such
maps into $\Pi_d$. The closure $\calx_d$ of this space in $\Pi_d$ is 
often referred to as the  
{\em Drinfeld's compactification} of the space of degree $d$ maps from 
$\CC P^1$ to $X$ and will be called the space of {\em quasimaps} 
(following \cite{FFKM}). It is a (generally speaking --- singular)
irreducible projective variety of complex dimension $\dim X+2d_1+...+2d_r$.

The flag manifold is a homogeneous space of the group $SL_{r+1}(\CC)$ and of
its maximal compact subgroup $SU_{r+1}$. The action of these groups on $X$
lifts naturally to the spaces $\Pi,\ \Pi_d,$ and $\calx_d$. 
In addition to
this action the spaces $\Pi_d$ and the subspaces $\calx_d$ carry the 
circle action induced by the rotation of $\CC P^1$ defined by 
$(x:y)\mapsto (x: e^{i\phi }y)$. Thus the product group $G=S^1\times SU_{r+1}$
and its complex version $G_{\CC}=\CC^* \times SL_{r+1}(\CC )$ act on
the quasimap spaces. We will see later that $\calx_d$ have 
$G$-equivariant desingularizations $\tilde{\calx}_d$.

We denote by $P=(P_1,...,P_r)$ the $G$-equivariant line bundles over 
$\calx_d$ 
obtained by pulling back the Hopf bundles over the complex projective factors 
of $\Pi_d$.

The $G$-equivariant holomorphic Euler characteristic of a  
$G$-equivariant holomorphic
vector bundle $V$ over a compact complex manifold $M$
provided with a holomorphic $G$-action is defined as the character of
the virtual representation of $G$ in the alternating sum of the cohomology
spaces:
\[ \chi_G (V):=\str_G (H^*(M, V) ):=\sum_k(-1)^k \tr_G (H^k (M, V) ) .\]
It can be expressed in cohomological terms using the equivariant version of 
the Riemann-Roch-Hirzebruch theorem:
\[ Ch (\x_G (V)) =\int_M Ch_G(V) Td_G(M), \]
where $Ch_G$ and $Td_G$ are the equivariant Chern character and Todd 
class, and $ch$ on the left hand side is the Chern character
from $K^*_G(pt)=K^*(BG)$ (canonically isomorphic to a suitable completion of 
the character ring $Repr (G)$) to $H^*_G(pt)=H^*(BG)$.

\medskip

We are interested in computing $G$-equivariant Euler characteristics of the 
line bundles $P^{ z}=P_1^{z_1}...P_r^{z_r}$ over equivariant 
desingularizations of the map spaces $\calx_d$. In fact (see
section $1$) the result does 
not depend on the choice of desingularization. 
We encode the answers by the following generating function:

\begin{equation} \calg (Q, z, q, \L) := 
\sum_{d} Q^d \chi_G (H^*(\tilde{\calx}_d, P^z))
\label{G} \end{equation} 

Here $q$ and $\L=(\L_0,...,\L_r | \L_0...\L_r =1)$ are multiplicative
coordinates on $S^1$ and on the maximal torus $T^r$ of $SU_{r+1}$ respectively,
and $Q=(Q_1,...,Q_r)$ are formal variables. The formula makes sense for
integer values of $z=(z_1,...,z_r)$ but can be extended to, say, complex
values of $z$ by interpreting the right hand side by means of the 
Riemann-Roch-Hirzebruch formula. 

\medskip

The {\em finite-difference Toda operator}
\footnote{We are thankful to P. Etingof for this definition.} 
\begin{equation} 
 \hat{H}_{Q,q} := q^{\p/\p t_0}+q^{\p/\p t_1}(1-e^{t_0-t_1})+...+
q^{\p/\p t_r}(1-e^{t_{r-1}-t_r}) \label{Toda} \end{equation}
is composed from the operators of multiplication by $Q_i=e^{t_{i-1}-t_i}$
and from the translation operators 
\[
  q^{\p/\p t_j}: t_i\mapsto t_i+ \delta_{ij} \ln q.
\]
We can now state the main result of the paper.

\medskip

{\bf Theorem 1.} {\em
In each of the two groups of variables $Q$ and
$Q'$ the function  
\[ G(Q,Q'):=\calg (Q, \frac{\ln Q'-\ln Q}{\ln q}, q, \L) \] 
is the eigen-function of the finite-difference 
Toda operator:
\[ \hat{H}_{Q',q} G\ =\ (\L_0+...+\L_r) G\ =\ \hat{H}_{Q,q^{-1}} G\ .\]
}

\medskip

{\bf Remarks.}
(1) Theorem $1$ 
together with some, rather general, factorization property of the generating 
function $\calg $ uniquely determine the function $\calg$. We will prove
this in Section $2$, and also present there some explicit examples. 

(2) Taking 
\[
 q=\exp (-\h) = 1-\h+\h^2/2+...
\]
and assigning the degrees $\deg \h =1, \deg Q_i=2$ we perform the degree 
expansion of the finite-difference Toda operator:
\[ \hat{H}= (r+1)-\h \sum \frac{\p}{\p t_i} +[\frac{\h^2}{2}
 \sum \frac{\p^2}{\p t_i^2} - \sum e^{t_{i-1}-t_i} ]+ ... .\]
In degrees $1$ and $2$ the expansion spits out the {\em momentum}
and the {\em Hamilton} operators of the quantum Toda lattice thus
explaining the name of $\hat{H}$. In fact the Toda operator $\hat{H}$ 
can be included into a complete set of commuting finite-difference 
operators (``conservation laws'')
in close analogy with the case of quantum and classical Toda systems. 
A construction of such operators in terms of quantum groups
is described in Section $5$.

(3) As we will see in Section $2$, the generating function $G$
is a common eigen-function of the commuting conservation
laws of the finite-difference Toda system. This result is a $K$-theoretic
counterpart (in the case of the series $A_r$) of the theorem by B. Kim 
\cite{BK} characterizing intersection theory in spaces of
holomorphic maps $\CC P^1 \to G/B$ in terms of quantum Toda lattices
(see Section $5.1$ for more details).

(4) A conjecture generalizing Theorem $1$ to the case of flag manifolds
$X=G/B$ of arbitrary semi-simple complex Lie groups $G$ and intertwining
K-theory on moduli spaces of stable maps with representation theory of
quantum groups is explained in Section $5$.     

(5) The results of the paper were completed in the Summer $98$ and
reported by the authors at a number of conferences and seminars.
In this version of the paper we decided to leave the material of Section $5$
in the form close to the preliminary text written in $1998$.
In particular, we did not try to match the quantum group description 
of finite-difference Toda lattices given in Section $5$
with the (apparently very similar) construction that has become 
standard since then due to the paper \cite{PE} 
by P. Etingof. Perhaps partly motivated by our conjectures, the 
paper places the finite-difference theory of Whittaker functions
on foundations much more solid than those available to us there years ago. 

(6) Initially the conjecture proved in this paper served as a motivation
for developing basics \cite{Gi2, YL} of {\em quantum K-theory}  
--- a K-theoretic counterpart of quantum cohomology theory. 
Moreover, one can heuristically interpret the generating function $\calg $
as an object of Floer-type (or semi-infinite) K-theory on the loop
space $LX$ equipped with the $S^1$-action defined by the rotation of loops.
According to \cite{Gi0} this heuristics, applied in cohomology theory,
suggests existence of a D-module structure (which in the case of 
flag manifolds is identified by Kim's theorem \cite{BK} 
with the quantum Toda system). 
The same arguments in K-theory lead to a finite-difference
counterpart of the D-module structure. While results of the present paper
conform with this philosophy, the role of finite-difference equations 
in the general structure of quantum K-theory remains uncertain 
(see Section 5 (d) in \cite{Gi1} for a few more details on this issue).

(7) The proof of Theorem $1$ is presented in Sections $1$ -- $4$.
In Section $1$ we discuss independence of the function $\calg $
on desingularizations. In Section $2$ we use desingularizations of
the quasimap spaces based on
moduli spaces of stable maps in order to obtain the factorization 
and other properties of the function $\calg $ mentioned in the Remarks $1,3$. 
In Section $3$ we describe the {\em hyperquot schemes} ---  
another equivariant desingularizations of the spaces $\calx_d$ ---
and compute their equivariant canonical class. This result plays a key role in 
our derivation of Theorem $1$ given in Section $4$.

\section{Rational desingularizations} 

A germ $(M,p)$ of a complex irreducible algebraic variety 
is called {\em rational} if for any desingularization 
$\pi : (M',p')\to (M,p)$ 
all higher direct images of the structure sheaf of $M'$ vanish in a
neighborphood of  $p$:
\begin{equation} (R^k\pi_* \calo_{M'})_p = 0 \ 
\text{for all} \  k>0. \label{rational} \end{equation}
Of course, $R^0\pi_*\calo_{M'} = \calo_M$. 

A non-singular germ is rational. This is a rephrasing of the famous
Grothendieck conjecture \cite{AG} proved by Hironaka \cite{HH}
on the basis of his resolution of singularity theorem
and saying that (\ref{rational}) 
holds true for any $f: M'\to M$ which is a proper  
birational isomorphism of non-singular spaces. 

Another easy consequence \cite{RH} 
of the Hironaka theorem is that the condition
(\ref{rational}) in the definition of rational singularities is automatically
satisfied for all desingularizations if it is satisfied for one of them. 
\footnote{This is proved by applying Leray spectral sequences to
the commutative square formed by three desingularizations 
$M'\to M,\ M''\to M$ and $M'''\to M$ dominating the first two: 
$M'''\to M', M'''\to M''$. We are thankful to R. Hartshorne for teaching us 
this subject.} In particular, it makes obvious the fact that
the product of a rational singularity with a non-singular space has rational 
singularities.

Let us call a {\em rational desingularization} of a singular space $N$
a proper birational isomorphism $M\to N$ where $M$ is allowed to have rational 
singularities at the worst. 
Consider a compact irreducible complex variety $N$ and
two birational desingularization $g_i: M_i\to N,\ i=1,2$. For any vector 
bundle $V$ on $N$ we have
$\chi (M_1; g_1^*V) = \chi (M_2, g_2^*V)$. Indeed, if $f_i:M\to M_i$ is
a common desingularization (so that  $f_1\circ g_1 = f_2\circ g_2$) then 
\[ \chi (M, f_i^*g_i^*V) = \chi (M_i, (R^*f_i)_*\calo_M \otimes g_i^*V) =
 \chi (M_i, \calo_{M_i} \otimes g_i^*V) = \chi (M_i, g_i^*V). \]
  
In our applications, we take on the role of $N$ the quasimap spaces
$\calx_d$ of the space of maps $\CC P^1\to X$. 
The holomorphic Euler characteristics we are
really interested in are those for bundles on $\calx_d$ pulled back to
the graph spaces (see Section $2$). 
For flag manifolds $X$ the graph spaces are orbifolds, and it is 
important for us that, due to a theorem by Viehweg \cite{EV}, 
{\em singularities of orbifolds are rational.} 
          
We need however the following equivariant refinement of the above 
independence argument:

\medskip

{\bf Proposition 1.} {\em Let a compact connected Lie group $G$ act by
holomorphic transformations on a compact irreducible complex projective 
variety $N$ and on two its equivariant projective rational 
desingularizations $f_i: M_i\to N,\ 
i=1,2$. Then for any $G$-equivariant holomorphic vector bundle $V$ on $N$
we have $\chi_G (M_1, f_1^*V) = \chi_G (M_2, f_2^*V)$.}

\medskip

{\em Proof.} The equivariant holomorphic Euler characteristic is a character
of a virtual representation of $G$ and is determined by its restriction to 
the maximal torus in $G$. Therefore we may assume that $G$ is such a torus.
Furthermore, the characters are determined by their (sufficiently high order)
jets at the point $1 \in G$. These jets have the following interpretation
in equivariant K-theory. 

Let $BG^{(n)}:=(\CC P^n)^s$ be the finite-dimensional approximations to the 
classifying space $BG=(\CC P^{\infty})^s$ of the $s$-dimensional torus $G$.
Let $\pi_i: (M_i)_G^{(n)} \to BG^{(n)}$ be $M_i$-bundles associated
with the restriction of the universal $G$-bundle to $BG^{(n)}$. Similarly,
consider the associated $N$-bundle $\pi: N^{(n)}_G\to BG^{(n)}$,
the bundle maps $F_i: (M_i)_G^{(n)} \to N_G^{(n)}$ induced by the 
equivariant maps $f_i: M_i\to N$ and the vector bundle $V^{(n)}$ over
$N_G^{(n)}$ associated with the equivariant bundle $V$ over $N$.
Then the K-theoretic push-forwards $(\pi_i)_* (F_i^*V^{(n)})$
(which are elements in the Grothendieck group $K^*(BG^{(n)}$ defined 
as the alternated sums $\sum (-1)^k\sum R^k(\pi_i)_*(...)$ of higher direct
images) coincide with suitable jets of the characters $\chi_G (M_i, f_i^*V)$.
\footnote{The coincidence, tautological in topological equivariant K-theory 
\cite{AS}, holds true in the context of algebraic geometry for non-singular
complex manifolds due to the compatibility of 
algebraic-geometrical and topological K-theoretic push-forwards (see
\cite{BFM}). We need here a more general statement applicable
to possibly singular spaces $M$. We don't have a suitable general reference
and assume that $M$ are projective instead. Then one can use an equivariant
embedding of $M$ into a projective space $P$ in order to push-forward $f^*V$ 
from $K^0_G(M)$ to $K^0_G(P)$ and then apply the coincidence in question for
the non-singular space $P$.} 
Moreover, these elements are uniquely determined by their K-theoretic
Poincare pairing with elements of $K^*(BG^{(n)})$, i. e. by
holomorphic Euler characteristics of the form
\[ \chi ((M_i)_G^{(n)}, F_i^*(V^{(n)}\otimes \pi^* W))\ ,\]
where the vector bundles $W$ run a basis in $K^*(BG^{(n)})$.

The bundles $\pi_i $ are locally
trivial with fibers $M_i$ having only rational singularities.
Therefore their total spaces $(M_i)_G^{(n)}$ have only rational singularities
too. The proposition follows now from its non-equivariant version
applied to the rational desingularizations $F_i$ of the spaces $N_G^{(n)}$.

\section{Graph spaces and factorization}

{\bf 2.1.} 
The graph of a degree $d$ holomorphic map $\CC P^1 \to X$ is a genus $0$
compact holomorphic curve in $GX:= \CC P^1\times X$ of {\em bi-degree}
$(1,d)$, \this of degree $1$ in projection to $\CC P^1$ and of degree $d$
in projection to $X$. 
We define the {\em graph space} $GX_d$ as the 
moduli space of genus $0$, unmarked stable maps to $GX$ of bi-degree $(1,d)$.
For $X=G/B$, the graph spaces $GX_d$,
according to \cite{MK, BM, FP}, are compact complex projective
\footnote{See, for instance, \cite{FP} where projectivity of moduli spaces
of stable maps to complex projective manifolds is proved. We need this property
only to assure that Proposition $1$ applies to the graph spaces considered
as equivariant rational desingularizations of the quasimap spaces.} 
algebraic
orbifolds of dimension $2d_1+...+2d_r+\dim X$. They provide therefore 
compactification of spaces of degree $d$ holomorphic maps $\CC P^1\to X$
and inherit the action of $G_{\CC}=S^1_{\CC}\times SL_{r+1}(\CC)$ from the
componentwise action on target space.  

The natural birational isomorphism between the graph spaces and
quasimap spaces is actually defined by a regular map (see 
for example \cite{Gi1}) which can be described as follows. 
A bi-degree $(1,d)$
rational curve in $\CC P^1\times X$ projected to $\CC P^1\times \CC P^{n_i-1}$
by the Pl\"ucker map consists of the graph $\S_0$ of a degree $m_0\leq d_i$ 
map $\CC P^1\to \CC P^{n_i-1}$ and a few ``vertical'' curves $\S_j$ of 
bi-degrees $(0, m_j)$ with $\sum m_j=d_i-m_0$,  attached to the graph. 
The graph component is given by $n_i$ mutually prime binary forms of degree 
$m_0$. Denote $\z_j\in \CC P^1$ the images of the vertical curves $\S_j$ in 
projection to $\CC P^1$. Multiplying the binary forms by the common factor 
with roots of multiplicity $m_j$ at $\z_j$ we obtain the degree $d_i$ 
vector-valued binary form which specifies the image of our curve in 
$\CC P^{(d_i+1)n_i-1}$. The regular map
\[ \mu: \ GX_d \to \Pi_d \] 
is defined by the above construction applied to each component of the
Pl\"ucker embedding.

\medskip

{\bf 2.2.}
The following result allows to separate the variables $Q$ and $Q'$ in the
generating function $G (Q, Q')$. Let $p_1,...,p_r$ denote Hopf line bundles
on $\Pi$ pulled back to the flag manifold $X$ by the Pl\"ucker embedding,
and $p^z=p_1^{z_1}...p_r^{z_r}$ be their tensor product.  
Denote by $\lan V, W\ran := \chi_G (X; V\otimes W)$ 
the K-theoretic Poincar\'e pairing on $K_G^*(X)$.

\medskip

{\bf Proposition 2.} {\em  
\[ G (Q, Q')=\lan J(Q', q) p^{\ln Q'/\ln q}, p^{-\ln Q/\ln q}
J(Q, q^{-1}) \ran, \]
where $J$ is a suitable formal $Q$-series with 
coefficients in $K_G^*(X) \otimes \QQ (\L, q)$ 
{\em (described below) }.}

\medskip

In the description of the series $J$, and in the proof of Proposition $2$
we will encounter the following concepts standard in Gromov -- Witten
theory (see for instance \cite{MK, BM} for their definitions and properties).
We will use the symbol $X_{m,d}$ for the moduli space of genus $0$ degree $d$
stable maps $f: (\S, \e_1,...\e_m) \to X$ with $m$ marked points 
(and we intend to avoid the notation $(S^1\times X)_{0,(1,d)}$ 
for the graph spaces $GX_d$).
Let $L$ denote the {\em universal cotangent line} bundle over 
$X_{1,d}$ formed by the cotangent lines $T^*_{\e}\S$ to the curves $(\S,\e)$
at the marked point $\e$ (see \cite{Gi2} for a discussion of these line
bundles in the context of K-theory). Let $\ev_*: K_G^*(X_{1,d})\to K^*_G(X)$ be
the K-theoretic push-forvard by the map $\ev: X_{1,d}\to X$ defined by the
evaluation $f\mapsto f(\e)$ at the marked point. In these notations
\[  J(Q,q)=1+\frac{1}{1-q}\sum_{d\neq 0} Q^d\ev_* (\frac{1}{1-Lq}). \] 

\medskip

 {\em Proof.} It is based on localization to fixed points of the 
$S^1$-action on the graph spaces $GX_d$ and goes through with minor
modifications in the general setting of quantum K-theory described
in \cite{Gi2}.

Let $f: \S \to \CC P^1\times X$ be a bi-degree $(1,d)$ stable map which
represents in $GX_d$ a fixed point of the $S^1$-action. Then $f$ consists
of the graph of a constant map $\CC P^1\to X$ and of two stable maps
$f_{\pm}: (\S_{\pm}, \e_{\pm}) \to X_{\pm}$ of bi-degrees $(0, d^{\pm})$,
$d^{+}+d^{-}=d$, with one marked point $\e_{\pm}$ each, 
attached to the graph at
the points $f_{\pm}(\e_{\pm})$. Here $X_{\pm}$ are two copies of $X$, namely 
the slices of the product $\CC P^1\times X$ over the fixed points $(1:0)$ and
$(0:1)$ of the $S^1$-action on $\CC P^1$. In the extremal cases $d^{+}=0$
or (and) $d^{-}=0$ only one (none) of the maps $f_{\pm}$ is present. A 
fixed point component in $GX_d$ is therefore identified with the suborbifold
in the product $X_{1,d^{+}}\times X_{1,d^{-}}$ of moduli spaces of stable maps
to $X$ with one marked point, given by the diagonal constraint 
$\ev(f_{+})=\ev(f_{-})$ in $X\times X$ for evaluations at the marked  points. 

According to the construction of the map $\mu : GX_d\to \Pi_d$, the image 
$\mu [f]$ is represented by an $r$-tuple of {\em monomial} binary 
vector-forms with zeroes of
orders $d^{+}$ and $d^{-}$ at $(1:0)$ and $(0:1)$ respectively. 
 
The $S^1$-fixed points in $\Pi_d$ represented by these vector-forms
form a fixed point submanifolds isomorphic to $\Pi$. We denote it by  
$\Pi_d^{(d_{+})}$ as the fixed point set $\Pi_d^{S^1}$ consists of one 
copy of $\Pi$ for each $0\leq d_{+}\leq d$. The images $\mu[f]$ of the
above fixed from $GX_d$ form a copy of $X$ Pl\"ucker-embedded into 
$\Pi_d^{(d_{+})}$. 
\footnote{We will need this information and notations also in
Section $4$.} This implies that
the $S^1$-equivariant line bundles $P=(P_1,...,P_r)$ restricted to the fixed 
point component coincide with $\ev^*(p)\otimes q^{d^{+}}$. Here 
$p=(p_1,...,p_r)$ is the list of Hopf line bundles on $X$, and 
$q^m=(q^{m_1},...,q^{m_r})$ specify the $S^1$-actions on the $r$-tuple of
trivial line bundles.

Each of the vertical curves $f_{\pm}$ contributes a
$2$-dimensional summand into the conormal bundle to the fixed point component.
The two infinitesimal deformations of $[f_{\pm}]$ breaking the $S^1$-invariance
correspond to the shift of the vertical curve away from the slice $X_{\pm}$ and
to the smoothening of $\S $ at the nodal point $\e_{\pm}$. These deformations
contribute respectively the factors $(1-q^{\pm 1})$ and 
$(1-L_{\pm}\otimes q^{\pm 1})$ to the denominator of the Bott--Lefschetz 
localization formula. Here $L_{\pm}$ is the universal cotangent line bundle
over $X_{1, d^{\pm}}$ formed by cotangent lines $T^*_{\e_{\pm}}\S_{\pm}$.
The contribution of the fixed point component to the localization formula
can be therefore written as
\[ \lan (\ev_{+})_* [(1-q)^{-1}(1-L_{+}q)^{-1}]\ p^zq^{d^{+}z}, \ 
(\ev_{-})_* [(1-q^{-1})^{-1}(1-L_{-}q^{-1})^{-1}] \ran ,\]
(where however the factor $(\ev_{\pm})_*[...]$ should be omitted if 
$d^{\pm}=0$).
We use here transversality of $\ev_{+}\times \ev_{-}$ to the diagonal 
$i:\D\subset X\times X$, which allows us to replace the push-forward to
$X_{1,d^{+}}\times X_{1,d^{-}}$ of the structure sheaf of the fixed point 
component by $(\ev_{+}\times \ev_{-})^*(i_*(O_{\D}))$.

Summing the contributions over all $(d^{+},d^{-})$ with the weights 
$Q^{d^{+}+d^{-}}$ we find
\[ \calg = \lan J(Qq^z, q), p^z J(Q, q^{-1})\ran . \]
It remains to recall that $G(Q,Q')$ is transformed to $\calg $ by the 
substitution $Q'=Qq^z$.

\medskip

{\bf 2.3.} Proposition $2$ shows that Theorem $1$ has the 
following reformulation.

\medskip

{\bf Theorem 2.} {\em The $K^*_G(X)$-valued 
vector-series $ p^{\ln Q/\ln q} J(Q, q)$ is 
the eigen-vector of the finite-difference Toda operator $\hat{H}_{Q,q}$
with the eigen-value $\L_0^{-1}+...+\L_r^{-1}$.}
        
\medskip 

 The series $J$ turns into $1$ when reduced modulo $Q$. 
In particular,
for the manifold $X$ of complete flags in $\CC ^{r+1}$ application of
the operator $\hat{H}_{Q,q}$ to $p^{\ln Q/\ln q}J$
yields, modulo $Q$, the factor $p_1+p_1^{-1}p_2+...+p_{r-1}^{-1}p_r +p_r^{-1}$.
The factor represents the trivial bundle $\CC ^{r+1}$ and is thus equal to
$\L_0^{-1}+...+\L_r^{-1}$ in $K^*_G(X)$.

\medskip

 Forgetting the space $\CC^i, \ i=1,...,r$, in the flag 
$\CC^1\subset ...\subset \CC ^r\subset \CC^{r+1}$ defines $r$ projections
$X\to X^{(i)}$ to partial flag manifolds with the fiber $\CC P^1$. 
The degrees ${\bf 1}_1,...,{\bf 1}_r$ of the fibers considered as rational 
curves in $X$ form the basis in $H_2(X)$ dual to the basis 
$(-c_1(p_1),..., -c_1(p_r))$ in $H^2(X)$ (and are 
represented by the monomials $Q_1,...,Q_r$ in our generating series). 
This identifies 
$X^{(i)}$ with the moduli space $X_{0, {\bf 1}_i}$, the projection 
$X\to X^{(i)}$ --- with the forgetting map $\ft: X_{1, {\bf 1}_i}\to 
X_{0,{\bf 1}_i}$, and shows that the evaluation map
$\ev: X_{1,{\bf 1}_i} \to X$ is an isomorphism. This information
about curves of minimal degrees in $X$ allows us to compute $J$ modulo $(Q)^2$:
\[ J=1+\frac{1}{1-q}[\frac{Q_1}{1-p_1^2p_2^{-1}q}+...+
\frac{Q_i}{1-p_{i-1}^{-1}p_i^2p_{i+1}^{-1}q}+...+
\frac{Q_r}{1-p_{r-1}^{-1}p_r^2q}]\ + \ o(Q). \]
The finite-difference equation $\hat{H} I=(\sum \L_j^{-1}) I$ for 
$I=p^{\ln Q/\ln q} J$ is equivalent to the recursion relation
\[ [ p_1 (q^{d_1}-1) +... +p_ip_{i-1}^{-1}(q^{d_i-d_{i-1}}-1)+...
+p_r^{-1}(q^{-d_r}-1)] J_d = \]  
\begin{equation} p_2p_1^{-1} q^{d_2-d_1} J_{d-{\bf 1}_1}+...+
p_r^{-1}q^{-d_r} J_{d-{\bf 1}_r} \label{rec} \end{equation}
for the coefficients $J_d=(1-q)^{-1}\ev_*(1-qL)^{-1}$ of the series
$J=\sum J_d Q^d$. It is straightforward now to verify the relation for 
$d={\bf 1}_i$. 

\medskip

The relation (\ref{rec}) recursively determines the coefficients $J_d$ 
unambiguously, which implies 

\medskip

{\bf Corollary 1.} {\em The power $Q$-series $J$ is uniquely determined by
Theorem $2$ and by the constant term $J|_{Q=0}=1$.}

\medskip

Moreover, the uniqueness
argument applies to arbitrary eigen-functions of $\hat{H}$.
Let $\hat{D}$ be any finite-difference operator with 
coefficients polynomial in $Q$ which commutes with $\hat{H}$,
and let $I=p^{\ln Q/\ln q} \sum_{d\geq 0} I_d Q^d$ be a power 
$Q$-series with vector-coefficients $I_d\in K^*_G(X)$.   
The following statement is the finite-difference version of 
Kim's lemma \cite{BK} important in quantum cohomology theory of
flag manifolds.

\medskip

{\bf Lemma.} {\em If $I$ is an eigen-function of $\hat{H}$: 
$\hat{H} I = \L_{H} I$,  and $\hat{D} I $ is proportional to $I$ modulo $Q$:
$\hat{D} I \equiv \L_{D} I \mod Q$, 
then $I$ is an eigen-function of $\hat{D}$: $\hat{D} I =\L_{D} I$.}

\medskip

Indeed, $\hat{H} \hat{D} I = \hat{D} \hat{H} I =\L_{H} (\hat{D} I)$ 
and therefore
$\hat{D} I$ is an eigen-function of $\hat{H}$ with the same constant
term $\L_{D} I_0$ as $\L_{D} I$ and thus coincides with it due to the
uniqueness argument.   

\medskip

The conservation laws of finite-difference Toda systems discussed in
Section $5$ satisfy the hypotheses of the lemma, and we obtain

\medskip

{\bf Corollary 2.} {\em The vector-function $I$ (as well as the series
$G$) is a common eigen-function of the commuting conservation laws
of the finite-difference Toda system.}

\medskip 

{\bf 2.4. Example: $r=1$.} 
The problem of computing the series (\ref{G})
can be generalized to arbitrary compact K\"ahler manifold $X$. 
In the case $X=\CC P^r$ the quasimap spaces coincide with the projective spaces
$\CC P^r_d$ and are non-singular. The corresponding series (\ref{G}) can be
computed immediately by the holomorphic Bott -- Lefschetz formula. We have
$ \sum_{d=0}^{\infty} Q^d \chi_G (\CC P^r_d; P^{\otimes z})  =$
\[ =  -\sum_{d=0}^{\infty} \frac{Q^d}{2\pi i} 
\oint_{P\neq 0} \frac{P^{z-1} \ dP}
{\prod_{j=0}^{r} \prod_{m=0}^d (1-P\L_jq^{-m})}= 
\lan J (Q q^z, q), p^z J(Q, q^{-1}) \ran, \] 
where 
\[ \lan \Phi , \Psi \ran = -\frac{1}{2\pi i} \oint_{p\neq 0} 
\frac{\Phi (p) \Psi(p) dp}
{p \Pi_{j=0}^r (1-p\L_j)}, \ \ 
 \text{and}\ \ J=\sum_{d=0}^{\infty} \frac{Q^d}{\prod_{j=0}^r \prod_{m=1}^d
(1-p\L_jq^{m})} .\]
The contour of integration here includes all poles except $0$. The
Hopf line bundle $p$ satisfies the relation $(1-p\L_0)...(1-p\L_r)=0$
in $K^*_G(\CC P^r)$.

In the non-equivariant limit $\L_0=...=\L_r=1$ the 
vector-function $I:=p^{\ln Q/\ln q} J(Q, q)$ satisfies the 
finite-difference equation $D^{r+1} I = Q I$ (here $D\ I (Q):= I(Q)-I(qQ)$)
with the symbol resembling the famous relation in the quantum cohomology
algebra of $\CC P^r$.

In the special case of the manifold $\CC P^1$ of complete flags in $\CC ^2$
the series $I:=p^{\ln Q/\ln q} J  $ at 
$Q=\exp (t_0-t_1)$ reads
\[ I= p^{\frac{t_0-t_1}{\ln q}} \sum_{d=0}^{\infty} \frac{e^{d(t_0-t_1)}}
{\prod_{m=1}^d\ (1-p\L_0q^m)(1-p\L_0^{-1}q^m)} \]
and modulo $(1-p\L_0)(1-p\L_0^{-1})=0$ satisfies 
\[ [ q^{\p/\p t_0}+ q^{\p /\p t_1}(1-e^{t_0-t_1}) ] \ \ I = 
(\L_0+\L_0^{-1})\ I. \]

{\bf 2.5. Example: $r=2$.} The space of complete flags in $\CC ^3$ 
coincides with the incidence relation $(\text{line})\ \subset \ 
(\text{hyperplane})$ in $\CC P^2\times \CC P^{2*}$. In this case the series
$I:=p^{\ln Q/\ln q} J$ still can be written explicitly in terms of the algebra
$K_G^*(\CC P^2\times \CC P^{2*})$ described by the relations
\[ (1-p_1\L_0)(1-p_1\L_1)(1-p_1\L_2)=0,\ (1-\frac{p_2}{L_0})(1-\frac{p_2}{L_1})
(1-\frac{p_2}{L_2})=0,\ \L_0\L_1\L_2=1 .\]
We have
\[ I=p_1^{\frac{t_0-t_1}{\ln q}}p_2^{\frac{t_1-t_2}{\ln q}}
\sum_{d_1,d_2=0}^{\infty} \frac{e^{d_1(t_0-t_1)+d_2(t_1-t_2)}\ 
\prod_{m=0}^{d_1+d_2} (1-p_1p_2q^m)}
{\prod_{j=0}^2 [ \prod_{m=1}^{d_1} (1-p_1\L_jq^m) \prod_{m=1}^{d_2} 
(1-p_2\L_j^{-1}q^m) ]}, \]
It is not hard to check directly that modulo the relations the series 
satisfies
\[ \hat{H} I = (\L_0^{-1}+\L_1^{-1}+\L_2^{-1}) I, \ \text{where}\  
 \hat{H}=q^{\p/\p t_0}+q^{\p/\p t_1}(1-e^{t_0-t_1})+q^{\p/\p t_2} 
(1-e^{t_1-t_2}),\]
and that $ \calg = \lan I(Qq^z, q), I(Q, q^{-1}) \ran = $ 
\be \label{f3} 
 \sum_{d_1,d_2=0}^{\infty} \frac{Q_1^{d_1}Q_2^{d_2}}{(2\pi i)^2}
\oint_{P_1\neq 0}\oint_{P_2\neq 0}
\frac{P_1^{z_1-1} P_2^{z_2-1} \ \ \prod_{m=0}^{d_1+d_2} (1-P_1P_2q^{-m}) \ \ 
dP_1 \w dP_2 }
{\prod_{j=0}^2 [ \prod_{m=0}^{d_1} (1-P_1\L_jq^{-m})\ \prod_{m=0}^{d_2}
(1-P_2\L_j^{-1}q^{-m}) ]} .\end{equation}

\section{Hyperquot schemes and the canonical class}

{\bf 3.1.} A degree $d$ holomorphic map $\CC P^1\to X$ defines a flag of 
subbundles $E^1\subset ...\subset E^r\subset \CC^{r+1} $ in the trivial bundle
over $\CC P^1$ of degrees $c_1(E^1)=-d_1, ..., c_1(E^r)=-d_r$. The {\em 
hyperquot scheme} $HQ_d$ is defined as the moduli space of the diagrams 
$E^1\to ... \to E^r\to E^{r+1}=\CC^{r+1}$ where the morphisms 
$E^i\to E^{i+1}$ of vector bundles
are injective almost everywhere on $\CC P^1$. According to \cite{AB, CF} the
hyperquot schemes are compact non-singular algebraic manifolds. 

More generally, the 
construction applies to partial flag manifolds and to grassmannians (in which
case one obtains Grothendieck's quot-schemes studied in \cite{SS}) and thus 
provides some non-singular compactifications of spaces of parameterized 
rational holomorphic curves in these spaces.  In the case
of degree $-d$ line subbundles in $ \CC^{N}$ the quot-scheme coincides
with the projectivization $\CC P^N_d$  of the space 
of vector-valued binary forms described in the Introduction.

 Replacing the flag $E^1\to ...\to E^r\to \CC^{r+1}$ by the
top exterior powers $\w ^iE^i \to \w ^i\CC^{r+1}$ we obtain a natural map 
$ \l : HQ_d \to \Pi_d:=\Pi_{i=1}^r \CC P^{n_i-1}_{d_i} $
to the product of the quot-schemes. The image of this map is the quasimap
space $\calx_d$. According to \cite{AB,CF}, the 
hyperquot schemes are non-singular compact algebraic manifolds. 
Therefore the map 
\[ \l : HQ_d \to \calx_d \subset \Pi_d \] 
provides an equivariant desingularization of the quasimap space. 
\footnote{In fact \cite{AK} it is a small resolution.}

\medskip

{\bf 3.2.} Denote $K_d$ the canonical line bundle of the hyperquot scheme
$HQ_d$. The following description of $K_d$ in terms of the pull-backs
$P_1,...,P_r$ of the Hopf line bundles from $\Pi_d$ will play a key role 
in the proof of Theorem $1$.

\medskip

{\bf Theorem 3.} {\em The class of the canonical line bundle $K_d$
in $K^*_{G} (HQ_d)\otimes \QQ$ coincides with the pull-back by $\l $ of  
\[ q^{-k_d} P_1^{2+2d_1-d_2}P_2^{2-d_1+2d_2-d_3}...\ 
P_{r-1}^{2-d_{r-2}+2d_{r-1}-d_r}P_r^{2-d_{r-1}+2d_r}, \]
 where 
 \[ k_d = d_1+...+d_r+\sum \frac{(d_i-d_{i-1})^2}{2}. \]}

 {\em Proof.} The moduli spaces $HQ_d$ is naturally equipped with
the {\em universal flag} (see \cite{AB, CF}) 
$\E^1\subset ...\subset \E^{r+1}=\calo^{r+1}$ of locally free sheaves on
the product $\CC P^1\times HQ_d$ (such that restrictions to 
$\CC P^1\times \{ b \}$
are sheaves of sections of the bundles in the diagram 
$E^1\to ... \to E^{r+1}=\CC ^{r+1}$ 
representing $b\in HQ_d$). According to \cite{AB, CF}, the tangent sheaf 
$\calt_d $ of the hyperquot scheme  can be described  as the kernel
of the following surjection:
\[ \sum_i f_{i-1}^*\otimes \id - \id \otimes g_i :
\oplus_i Hom (\E^i, \calo^{r+1}/\E^i) \to \oplus_i 
Hom (\E^i, \calo^{r+1}/\E^{i+1}) .\]
Here $f_i: \E^i\to \E^{i+1}$ are the inclusions and $g_i: \calo^{r+1}/\E^i \to
\calo^{r+1}/\E^{i+1}$ are corresponding projections. 
In other words, the class of the
tangent bundle $\calt_d$ to $HQ_d$ in the Grothendieck group 
$K^*_{S^1\times G}(HQ_d)$ equals
\[ \oplus_{i=1}^r [\ Ext^0 (\CC P^1; \E^i, \E^{r+1}/\E^i) \ominus 
Ext^0 (\CC P^1; \E^i, \E^{r+1}/\E^{i+1}) \ ] .\]
It is easy to see from long exact sequences generated by
$0 \to \E^j \to \E^{r+1} \to \E^{r+1}/\E^j \to 0$  that 
$Ext^1 (\CC P^1; \E^j, \E^{r+1}/\E^j)=0$.  Thus the class of tangent bundle is 
represented by the $K$-theoretic push-forward along the projection 
$pr: \CC P^1\times HQ_d \to HQ_d$:  
\[ \calt_d = pr_* [ \ \oplus_i (\E^i)^* \otimes (\E^{i+1}-\E^i) \ ] . \]

We intend to compute the equivariant $1$-st Chern class of $\calt_d$ by
means of the relative Riemann -- Roch theorem:
\[ Ch ( pr_* V)= pr_* [ Ch (V) Td (\CC P^1)  ] ,\]
where $Ch$ and $Td $ are equivariant Chern character and Todd class 
respectively. Since fibers of $pr$ have dimension $1$, it suffices to
confine only terms of degree $\leq 2$ in $Ch (V)$ and $Td (\CC P^1)$ in order
to compute the $1$-st Chern class of $pr_* (V)$. Our computation goes through
due to the following ``miracle'': for a virtual bundle $V$ of the form
$V=\sum (E^i)^*(E^{i+1}-E^i)$, 
where $E^i$  have dimensions $i$ and $E^{r+1}$ is trivial, the degree 
$\leq 2$ terms of $Ch (V)$ depend only on $1$-st Chern classes $c_1(E^i)$.
Namely,
\[ Ch (V)=\frac{r(r+1)}{2}-2\sum c_1(E^i) +\frac{1}{2} 
\sum (c_1(E^{i+1})-c_1(E^i))^2 + ... \]
In our situation $E^{r+1}=\CC^{r+1}$ is topologically trivial, but carries
a non-trivial action of $SL_{r+1}(\CC)$. Yet $c_1(E^{r+1})=0$, but 
we get an extra summand $r c_2 (E^{r+1})$. This summand however is
a constant $const \in H^2(BSU_{r+1})$ and will disappear from our formulas
after integration over $\CC P^1$.  
 
Notice that $c_1(E^i)=c_1(\w^i E^i)$, and that the subsheaves 
$\w^i \E^i \subset \w^i \calo^{r+1}$ are exactly those which define the map 
$HQ_d \to \Pi_d$. This allows us to perform our computation in the product
$\CC P^1\times \Pi_d$ of projective spaces instead of $\CC P^1\times HQ_d$.

The equivariant cohomology algebra $H^*_{S^1}(\CC P^1)$ with respect to our
usual $S^1$-action is isomorphic to $\ZZ [\rho ,\h ]/(\rho (\rho +\h))$
where $\h =-\ln q$ is the generator of $H^*_{S^1}(pt)=H^*(BS^1)$, 
and $-\rho $ is the
equivariant $1$-st Chern class of the Hopf line bundle over $\CC P^1$.
The cohomological push-forward $pr_*: H^*_{S^1}(\CC P^1 ) \to H^*_{S^1}(pt) $ 
is given by $pr_* (1)=0$, $pr_* (\rho )=1$. Computing the equivariant 
$1$-st Chern class $c=2\rho +\h$ we expand the equivariant Todd class:
\[ Td (\CC P^1) =1+\frac{c}{2}+\frac{c^2}{12}+... = 
1 + (\rho -\frac{\h}{2}) + \frac{\h^2}{12} +...   \]

Consider now the quot-scheme
$ {\mathbb P}_d:=Proj ( Ext^0(\CC P^1; \calo (-d), \calo^N ) )$ 
corresponding to
degree $d$ maps $\CC P^1 \to Proj (\CC ^N)$. The equivariant cohomology algebra
$H^*_{S^1\times G} ({\mathbb P}_d)$ (with respect to any linear $G$-action on 
$\CC ^N$) is generated over the coefficient ring $H^*(BS^1\times BG)$ 
by the equivariant $1$-st Chern class $H=c_1(P)$ of the Hopf line bundle $P$. 
The tautological composition 
$\calo (-d) \otimes Hom ( \calo (-d), \calo^N ) \to \calo^N$ defines 
the universal rank $1$ subsheaf $\w \subset \calo ^N$ on 
$\CC P^1 \times {\mathbb P}_d$ where therefore $\w =\calo_{\CC P^1}(-d)\otimes P$.
Thus $c_1(\w )=H-d\rho $.

Now we apply the above construction to each factor in $\Pi_{d_1,...,d_r}$ by 
putting $\CC^N=\w^i \CC^{r+1}$, $d=d_i$, $i=1,...,r$, and conclude that 
the equivariant $1$-st Chern classes $c_1(\w^iE^i)$ are represented in
$\CC P^1\times \Pi_d$ by $H_i-d_i\rho $ where $H_i=c_1(P_i)$. We compute
the degree $2$ term in $Ch (V) Td (\CC P^1)$:
\[ \frac{r(r+1)}{24}\h^2-(2\rho+\h)\sum (H_i-d_i\rho)+
\frac{1}{2}\sum (H_{i+1}-H_i-(d_{i+1}-d_i)\rho )^2 + const .\]
Replacing powers of $\rho $ in this formula by their push-forwards 
$pr_*(1)=0$, $pr_*(\rho )=1$, $pr_*(\rho^2)=pr_*(-\h \rho)=-\h$,
$pr_*(const)=0$, we find
\[ c_1(\calt_d)=-\h [\sum d_i +\frac{1}{2}\sum (d_{i+1}-d_i)^2]  
-2\sum H_i -\sum (H_{i+1}-H_i)(d_{i+1}-d_i) .\] 
Since $c_1(K_d)=-c_1(\calt_d)$, $\h =-c_1(q)$ and $H_i=c_1(P_i)$, 
we finally conclude that (at least over $\QQ $)
\[ K_d = q^{-k_d} \prod_i P_i^2 (P_{i+1}P_i^{-1})^{d_{i+1}-d_i} =
q^{-k_d} \prod_i P_i^{2-d_{i-1}+2d_i-d_{i+1}}. \]  

\medskip

{\bf 3.3.} Another useful property of the hyperquot schemes is
that {\em fixed points of the action on $HQ_d$ of the maximal torus 
$S^1\times T \subset G$ are isolated}. More precisely, 
a fixed point of the torus $S^1\times T$ action on $HQ_d$. Such a 
fixed point is uniquely determined by the following data:

(i) A permutation $\s \in S_{n+1}$ specifying a fixed point of the torus $T$
action on the flag manifold $X$.

(ii) A pair $\D_{+}, \D_{-}$ of lower-triangular matrices with non-negative
integer entries $m_{ij}$, $1\leq j\leq i\leq r$ satisfying
\[ 0\leq m_{i1}\leq ...\leq m_{ii},\ i=1,...,r,\ \text{and}    
\ \sum_{i=j}^r (m^{+}_{ij}+m^{-}_{ij}) =d_{r+1-j},\ j=1,...,r. \]

Indeed, let the flag of subsheaves 
$\E^1\subset ...\subset \E^{r+1}=\calo^{r+1}$ on $\CC P^1$
represent a fixed point of the torus $S^1\times T$ action on $HQ_d$. 
At generic points of $\CC P^1$ the flag of subsheaves determines a 
flag of subspaces in $\CC^{r+1}$ which may not vary along $\CC P^1$ 
($S^1$-invariance) and thus coincides with one of $(r+1)!$ coordinate flags
in $\CC^{r+1}$ ($T$-invariance). Let $(e_0,...,e_r)$ be the standard basis
in $\CC ^{r+1}$. Consider for example the $T$-invariant flag formed by the 
coordinate subspaces $\CC^{r+1-j}:=Span (e_j,...,e_r)$ (all other $T$-invariant
flags are obtained by permutations $(e_{\s(0)},...,e_{\s(r)})$ of the basis).
Outside the fixed point set $(1:0),(0:1)$ of $S^1$-action on $\CC P^1$
the sheaf $E^{r+1-j}$ coincides with the sheaf of vector-functions
$ \calo e_j\oplus \calo e_{j+1}\oplus ... \oplus \calo e_r $ 
($S^1$-invariance).
It remains to describe the $S^1\times T$-invariant flags of subsheaves near 
the fixed points. It is easy to see that such a flag is equivalent to one
of the following, described by the matrices $\D_{+}$ and $\D_{-}$ near 
$(1:0)$ and $(0:1)$ respectively. Let $\z $ be the coordinate on $\CC P^1$ near
$(1:0)$, and $(\z^m)$ denote the ideal generated by $\z^m$ in the local algebra
of functions on $\CC P^1$ near at this point. In a neighborhood $U$ of 
$\z=0$ put
\[ \E^{r+1-j} |_U:=\oplus_{i=j}^r (\z^{m_{ij}}) e_i \ . \]
The conditions on the matrices $\D_{\pm}$ guarantee the subsheaves form a flag 
and that their degrees are equal to $d_{r+1-j}$. 
  
\medskip

{\bf Remark.} One can continue along these lines and write down explicitly
Bott -- Lefschetz fixed point localization formulas on $HQ_d$. In particular, 
one can easily recover this way the factorization property of the function $G$
described by Proposition $2$. On the other hand, we were unable to see
how combinatorics of the localization formulas  
(which in principle determine the function $G$) implies the recursion
relation (\ref{rec}). In the proof of Theorem $1$ given in the
next section we choose a different way which refers to
localization formulas in the projective spaces $\Pi_d$ instead. 
What we need to know about the hyperquot schemes is only Theorem $3$ and the 
very fact that fixed points in $HQ_d$ are isolated.  
   
\section{Localization and recursion}

{\bf 4.1.} The recursion relation (\ref{rec}) can be restated as an identity
between some elements in $K^*_{S^1\times G}(\Pi_d)$. Let 
\[ O_d:=\mu_*(\calo_{GX_d}) \]
be the $K$-theoretic push-forward of the trivial line bundle over the graph 
space to the product $\Pi_d$ of projective spaces along the map 
$\mu : GX_d\to \Pi_d$ described in Section $2.1$.
The elements $O_d$ with different $d$ can be compared to each other via
the inclusions 
\[ \phi^{(i)}: \Pi_{d-{\bf 1}_i} \subset \Pi_d \]
defined in terms of $r$-tuples 
$(f_1,...,f_r)$ of vector-valued binary forms by the formula 
$f_j(x,y)\mapsto f_j(x,y)y^{\d_{ij}}$. As it is easy to see, say, from
localization to fixed points of $S^1$-action on $\Pi_d$, the pull-back by
$\phi^{(i)}$ transforms $P_j$ to $q^{\d_{ij}} P_j$. Let 
\[ O_d^{(i)}:= \phi^{(i)}_* O_{d-\1_i},\ i=1,...,r, \] 
be the K-theoretic push-forward of $O_{d-{\bf 1}_i}$ by the inclusions. 
One can think of $O_d$ and $O_d^{(i)}$ as Laurent polynomials
of the generators 
$(P_1,...,P_r, \L_0,...,\L_r, q)$ in $K^*_{ G}(\Pi_d)$ 
defined modulo relations among them. Remembering the notation
$Q'_i=Q_i q^{z_i}$, we get $\hat{H}_{Q',q} G (Q, Q') = \sum_d Q^d \calh_d $
where 
\[ \calh_d\  =\  \chi_G (\Pi_d ; O_d
(P_1+P_2P_1^{-1}+...+P_r^{-1}) P^z ) 
 - q^{z_1}\chi_G (\Pi_{d-\1_1} ; O_{d-\1_1}P_1 P^z) \]
\[  - q^{z_2}\chi_G (\Pi_{d-\1_2}; O_{d-\1_2}P_2P_1^{-1} P^z)-...- 
q^{z_r}\chi_G (\Pi_{d-\1_r}; O_{d-\1_r}P_r^{-1} P^z) .\]
This shows that Theorem $1$ is equivalent to the following 
sequence of relations in $K^*_G(\Pi_d)$:
\begin{equation} H_d:= P_1 O_d + P_2P_1^{-1}(O_d-O_d^{(1)})+...+
P_r^{-1}(O_d-O_d^{(r)}) = (\sum \L_j^{-1}) O_d . \label{id} \end{equation}
Moreover, due to the results of Section $2$, the relations (\ref{id})
for all $d\leq d_0$ are equivalent to the recursion 
relations (\ref{rec}) for coefficients $J_d$ of the series $J$ for all
$d\leq d_0$. We are going to prove (\ref{id}) by induction on 
$|d|:=d_1+...+d_r$. 

As we know from Section $2.3$ the relation (\ref{id}) is true
for $|d| \leq 1 $.   

Let us now assume that (\ref{id}) is true for all $d'\neq d$ such that
$d'\leq d $ (componentwise). 
We compute the fixed point localizations of $H_d$ at the 
fixed point components $\Pi_d^{(d^+)}$ of the $S^1$-action.
More precisely, we call here 
the {\em localization} of a torus-equivariant class 
$A\in K_T^*(Y)$ the class $a$ in the (localized) equivariant $K$-ring of the
fixed point set $Y^T$ such that $j_*a=A$ under the embedding $j: Y^T\to Y$.
The localization is therefore characterized by the property 
$\chi_T(Y^T; a j^*B) = \chi_T (Y; A B)$ for any $B\in K^*_T(Y)$. 

Let us recall from Section $2.2$ that $S^1$-fixed points in $\Pi_d$ 
are represented by  
vectors of binary forms with {\em monomial} components proportional to 
$ x^{d_i^{-}}y^{d_i^{+}}$ (where $d^{-}=d-d^{+}$), 
that the set $\Pi_d^{(d^+)}$ formed by these fixed points is a copy of
$\Pi $ and that the restriction of $P_i$ to $\Pi_d^{(d^+)}$ equals 
$p_iq^{d_i^+}$. Due to Proposition $1$ we have,  
$O_d := \l_* (\calo_{HQ_d})=\mu_*(\calo_{GX_d})$. Due to the 
results of Section $2$ the localization
of $O_d$ to the fixed component belongs therefore to the image of 
$i_*: K_G^*(X) \subset K_G^*(\Pi)$ under the Pl\"ucker embedding and
is described in terms of $K^*_G(X)$ by the coefficients of the series
$J(Q,q)=\sum J_d(q)Q^d$ as $i_* [ J_{d^+}(q) J_{d^-}(q^{-1}) ]$.

Similarly, localizations of $O_d^{(i)}:=\phi^{(i)}_*O_{d-\1_i}$ coincide
with $i_*[ J_{d^+}(q) J_{d^{-}}(q^{-1})]$ since $\phi^{(i)}$ maps 
isomorphically the fixed point component $\Pi_{d-\1_i}^{(d^+-\1_i)}$ onto
$\Pi_d^{(d^{+})}$ for all $d^+\geq \1_i$. Combining, we find that
localizations of $H_d$ have the form $i_*[C_{d_+}(q) J_{d^-}(q^{-1})]$
where $C_d =$
\[ ( p_1 q^{d_1} +... +p_ip_{i-1}^{-1}q^{d_i-d_{i-1}}+...
+p_r^{-1}q^{-d_r}) J_{d} -  
p_2p_1^{-1} q^{d_2-d_1} J_{d-{\bf 1}_1}-...-
p_r^{-1}q^{-d_r} J_{d-{\bf 1}_r} .\]
Comparing with the recursion relation (\ref{rec}) we find that
the induction hypothesis implies vanishing of localizations of 
$H_d-(\sum \L_j^{-1}) O_d$
at all fixed point component $\Pi_d^{(d^+)}$ with $d^+  \neq d$
(and the remaining vanishing condition for $d^+=d$ coincides with 
(\ref{rec}) ).

{\bf Remark.}
Vanishing of the localization at the last fixed set $\Pi_d^{(d)}$ will
be derived with the use of the following polynomiality property. 
The classes $H_d$ and $O_d$ 
are defined in $K^*_G(\Pi_d)$ before of fixed point
localization (while their expression in terms of $J_{d^{\pm}}$ makes 
explicit use of it) and therefore     
\[ \chi_G (\Pi_d; P^z [H_d-(\sum \L_j^{-1})O_d])\ \in Repr (G) \]
is the character of a virtual representation of $G$ and is therefore
represented by a {\em Laurent polynomial}, that is a regular function
on the maximal torus of the group $G_{\CC}=S^1_{\CC}\times SL_{r+1}(\CC)$. 
When expressed in terms of the localizations $J_{d^+}(q)J_{d^{-}}(q^{-1})$ 
(which typically have lots of other poles in $q$ besides $q=0,\infty$) 
the polynomiality property yields serious constraints on $J_d$. Yet
it turns out (although we are not going to describe the details here) that 
the constraints are not powerful enough in order to uniquely 
determine all $J_d$, and we need additional geometrical
information about them.

\medskip

{\bf 4.2.}Consider further localizations $J_d^{\s}$ of the coefficient 
$J_d(q)\in K^*_G(X)$ to the $(r+1)!$ fixed points $\s\in X$ of the maximal 
torus $T\subset SU_{r+1}(\CC)$. The restriction are rational functions of 
$q$ and $\L$ and can be written in the form
\[   J_d^{\s} = \frac{ R_d^{\s}(q) }{S_d^{\s}(q)}, \]
where $R_d^{\s}, S_d^{\s}\in \QQ(\L ) [q]$ are polynomials in $q$ not vanishing
at $q=0$ simultaneously. We claim that in fact
\be \deg S^{\s}_d - \deg R^{\s}_d \geq k_d =
d_1+...+d_r+\sum \frac{(d_i-d_{i-1})^2}{2} \label{mgn} \end{equation}
(where we put $d_i=0$ for $i\leq 0$ and $i>r$).
This follows from general structure of Bott -- Lefschetz fixed point 
localization formulas applied to the hyperquot schemes and from Theorem $3$.

Indeed, consider the (isolated!) 
fixed points $(\s,\D)$ (see Section $3.3$) in $HQ_d$
mapped to the fixed point $\s \in \Pi_{d}^{(d)}$.
As we found in Section $4.3$, the coefficient $J_d^{\s}$
represents the localization of $O_d=\mu_*(\calo_{GM_d})$ at the $(r+1)!$
fixed points in $\Pi_d^{(d)}$. On the other hand we have 
$O_d=\l_*(\calo_{QH_d})$ due to Proposition $1$. Therefore the localization
coefficient $J_d^{\s}$ is equal to the sum 
\[ J_d^{\s} = \sum_{\D } \frac{1}{BL_{\s,\D}} \]
of the similar localization coefficients for $\calo_{HQ_d}$ at the 
fixed points $(\s,\D)$ mapped to $\s$. 
The Bott -- Lefschetz denominator $BL_{\s,\D}$ here
is the character of the torus $S^1\times T$-action on the exterior algebra
$\w^*\calt^*_{\s,\D}$ of the cotangent space to $HQ_d$ at the fixed point.
Therefore
\[ BL_{\s,\D} = \Pi_{\a} (1-q^{\v_{\a}}\L^{\x_{\a}}) .\]
where $(\v_{\a},\x_{\a})$ specify the characters of respectively $S^1$ and $T$
on the $1$-dimensional invariant subspaces in $\calt^*_{\s,\D}$ indexed by 
$\a$. 

Each fraction $1/(1-q^{\v}\L^{\x})$ with $\v >0$ has order $0$ at $q=0$ and
order $\v $ at $q=\infty $. Each fraction with $\v<0$ can be rewritten as
$-q^{-\v}\L^{-\x}/(1-q^{-\v}\L^{-\x})$ and has order $0$ at $q=\infty$ and
order $-\v$ at $q=0$. The product $BL_{\s,\D}^{-1}$ of such fractions
is uniquely written as a rational functions in $q$ of the form 
$Const \ q^M / S(q)$, where $M=-\sum \v_{\a}$ over all negative $\v_{\a}$,
$S$ is a polynomial in $q$, $S(0)=1$, and the degree
$\deg S = \sum_{\a} \v_{\a} + 2M$ of the denominator exceeds the degree $M$
of the numerator by $\sum_{\a} \v_{\a}$ at least. Clearing denominators
in the sum $J_d^{\s}$ of such rational functions we represent $J_d^{\s}$ as
a ratio $R_d^{\s}(q)/S_d^{\s}(q)$ of two polynomials in $q$ where 
\[ S_d^{\s}(0)=1\ \text{and}\  
\deg S_d^{\s}-\deg R_d^{\s} \geq \sum_{\a} \v_{\a} .\]

Notice that the sum $\sum_{\a} \v_{\a}$ coincides with the character of the 
$S^1$-action on
the top exterior power of $\calt^*_{\s,\D}$. Using Theorem $3$ together
with the fact that restrictions of $P_i$ to $\Pi_d^{(d)}$ are equal
to $q^{d_i}p_i$ we find  
\[ \sum_{\a} \v_{\a} = -k_d+\sum d_i(2-d_{i-1}+2d_i-d_{i+1}) = -k_d+2k_d=k_d.\]
regardless of the index $\D$. 

\medskip

{\bf 4.3.} 
With the estimate (\ref{mgn}) at hands we now complete the induction 
step as follows. The estimate shows that the total order of poles of 
$J_d^{\s}$ at $q\neq 0, \infty$ is $k_d$ at least. 

On the other hand, we claim that the   
localizations 
$H_d^{\s}-(\sum \L_j^{-1}) O_d^{\s}$ of \linebreak
 $H_d-(\sum \L_j^{-1})O_d$ at the fixed points of $S^1\times T^r$
with $d^+=d$ and any $\s\in X^T$ have no poles at $q\neq 0,\infty$.

Indeed, for any $z\in \ZZ^r$ the $G$-equivariant holomorphic Euler
characteristic of the sheaf $P^z [H_d-(\sum \L_j^{-1})O_d ]$
is a Laurent polynomial of $q$ and $\L$. By the induction
hypothesis localizations of this sheaf to all fixed points of $S^1$-action with
$d^+\neq d$ vanish, and thus the Euler characteristics is equal to
\[ q^{zd} \sum_{\s} \L_{\s}^z [ H_d^{\s}-(\sum \L_j^{-1})O_d^{\s}], \] 
The restrictions $\L_{\s}^z=\L_{\s (0)}^{-z_1}\L_{\s (1)}^{z_1-z_2}...
\L_{\s (r)}^{z_r}$ of $p^z$ to the $(r+1)!$ fixed points $\s$ (identified in
this formula with corresponding permutations) are linearly independent
exponential functions of $z$. This shows for each $\s$
\be H_d^{\s} -(\sum \L_j^{-1}) O_d^{\s} = \sum_{i=0}^r 
\L_{\s (i)}^{-1} q^{d_{i+1}-d_i}(J_d^{\s}-J_{d-{\bf 1}_i}^{\s})
-(\sum \L_i^{-1}) J_d^{\s}, \label{hs} \end{equation}
(where we put $C_{d-{\bf 1}_0}=0$) must be Laurent polynomials.

Using the estimate (\ref{mgn}) and the identities
\[ k_{d-{\bf 1}_i}-(d_{i+1}-d_i)=k_d-(d_i-d_{i-1}), \]
we see that (\ref{hs}) multiplied by $q^{\max_i (d_i-d_{i+1})}$ 
is a rational function of the 
form $R(q)/S(q)$ with $S(0)\neq 0$ with 
\[ \deg S-\deg R\geq k_d-\max (d_i-d_{i-1})-\max (d_i-d_{i+1}) .\]
This conclusion together with the property of (\ref{hs}) to have no finite
non-zero poles leads to a contradiction 
unless $R=0$ or $k_d\leq \max (d_i-d_{i-1})+\max (d_i-d_{i+1})$.
The inequality implies $|d|\leq 1$ (and turns into equality
for $d=0$ and $\1_i$.)
Indeed, denote the maximum and the minimum by $\a$ and $\b$. We have
\[ 0\leq \a+\b-k_d= \a+\b -\sum \frac{(d_i-d_{i+1})^2}{2}-|d| \leq 
\sum \a-\frac{\a^2}{2}+\b-\frac{\b^2}{2}-|d| \leq 1-|d|.\] 
Thus for $d_1+...+d_r> 1$ we are left with the only option $R=0$. This
completes the induction step and the proof of Theorem $1$. 

\medskip

{\bf Remarks.}  (1) In cohomology theory, the arguments parallel to
those used in Section $4.1$ 
(\this the general factorization and polynomiality properties plus
explicit description of spaces of curves of degrees $|d|\leq 1$)
would have been sufficient in order to determine the counterpart of
the series $J$ unambiguously and thus prove Kim's theorem \cite{BK}
mentioned in the Remark (3) in the Introduction (see also Section $5.1$ 
below). This follows from dimension
counting: for $|d|>1$ the dimension of $GX_d$ exceeds the degree
of the cohomology class analogous to $H_d-(\sum \L_j^{-1}) O_d$ by a margin
large enough in order to provide the necessary number of
additional linear dependencies among fixed point localizations.    

(2)
In $K$-theory, the dimensional argument is not available. Let us look however
at the integral formula (\ref{f3}) (yet conjecturally) 
representing $\calg $ for 
$r=2$. For $d_1, d_2$ large enough and $z_1,z_2 >0$ small enough the integrand 
has no poles with $P_1,P_2=0$ or $\infty$, and thus the integral is equal 
to $0$. If known 
{\em a priori}, this property would provide enough linear dependencies
between localizations of $H_d$ in order to complete our proof.

(3) The same property of the countur integral (see Section $2.4$) 
representing the Bott -- Lefschetz
formula for complex projective spaces follows {\em a priori} from 
the Kodaira -- Nakano vanishing theorem. It would be interesting to figure
out what kind of vanishing theorems would guarantee this property in
the case of the graph spaces $GX_d$.

(4) In the above proof we exploited another property which is manifest
in (\ref{f3}): the difference $d_1^2/2 + (d_1-d_2)^2/2+d_2^2/2+d_1+d_2$ 
between the degrees of the numerator and the denominator in the integrand
considered as a function of $q^{-1}$. The margin is explained
by Theorem $3$, and the whole argument resembles the famous proof \cite{AH}
of Atiyah -- Hirzebruch rigidity theorem of arithmetical genus
refined by the estimate of the equivariant canonical class. 

\section{Generalization to flag manifolds $G/B$}

An arbitrary complex semi-simple Lie group $G$ will
replace here the group $SL_{r+1}(\CC)$ of the previous sections.

\medskip

{\bf 5.1.} In order to generalize Theorems $1$ and $2$
 to flag manifolds 
$G/B$ of semi-simple complex Lie algebras $\Gg $ it is useful to recall 
corresponding results of quantum cohomology theory due mainly to B. Kim 
\cite{BK}. According to the Borel -- Weil construction, fundamental 
representations $V_1,..., V_r$ of the Lie algebra $\Gg $ of rank $r$ can
be realized in the spaces of holomorphic sections of suitable line bundles
over $G/B$ with the $1$-st Chern classes $h_1,...,h_r$. 
The sections define the {\em Pl\"ucker embedding} 
\[ X:=G/B \to \Pi := \Pi_{i=1}^r Proj (V_i^*), \]
which allows one to generalize the  construction
of the maps $\mu$ graph spaces from $GX_d$ to the products 
$\Pi_d$ of projective spaces.
The cohomological counterpart of the series (\ref{G}) is
\[ \calg^{H}=\sum_{d=(d_1,...,d_r)} Q_1^{d_1}...Q_r^{d_r} 
\int_{GX_d} e^{H_1z_1+...+H_rz_r}, \]
where $H_i$ are $S^1\times G$-equivariant $1$-st Chern classes of the
hyperplane line bundles on the product $\Pi_d$ of projective spaces pulled
back to the graph space $GX_d$ by the map $\mu $. The power $(Q,z)$-series 
$\calg_H$ with coefficients in $H^*(BS^1\times BG, \QQ)$ can be factored as
$( I_H(Qe^{\h z}, -\h^{-1}), I_H(Q, \h^{-1}) )$ where $( .\, , . )$ 
is the $G$-equivariant Poincar{\'e} pairing on $H^*_G(X)$, 
$I_H(Q,\h^{-1})=e^{(h\ln Q)/\h)}[\sum_d c_d Q^d]$ is a suitable series with 
coefficients $c_d\in H^*_G(X, \QQ (\h))$, and $\h $ is the generator of
$H^*(BS^1)$.
  
The main theorem in \cite{BK} says that {\em $I_H$ is a common 
eigen-function of $r$ commuting differential operators which form a complete 
set of quantum conservation laws of the quantum Toda lattice corresponding to
the Langlands-dual Lie algebra $\Gg '$}. 

For example, the Hamiltonian operator 
of the Toda lattice in the (self-dual) $sl_{r+1}$-case has the form 
\[ H=\frac{\h^2}{2} \sum_{i=0}^r \frac{\p^2 }{\p t_i^2} 
-\sum_{i=1}^r e^{t_{i-1}-t_i}, \]
and $H\ I =\frac{1}{2}(\l_0^2+...+\l_r^2)\ I$ (we identify here $H^*(BG)$ 
with the algebra of symmetric functions in $\l_0,...,\l_r,\ \sum \l_i=0$). 
As we mentioned in the Remark (2) in the Introduction, 
the equation for $I_H$ can be extracted from the finite-difference equation 
$\hat{H} I =(\sum \L_i^{-1}) I$ as the degree $2$ part in the following
approximation: put $q=e^{-\h }=1-\h+...$, $\L_i=e^{\l_i}=1+\l_i+...$ and 
use the grading $\deg \h =1, \deg \l_i=1, \deg Q_i=2 $. 

\medskip

In $K$-theory, we introduce
\[ \calg_K:= \sum_{d=(d_1,...,d_r)} Q_1^{d_1}...Q_r^{d_r} 
\chi_{S^1\times G} (GX_d; P_1^{z_1}...P_r^{z_r})) ,\]
where $P_i$ denote the Hopf line bundles on $\Pi_d$ pulled back to $GX_d$
by $\mu $. Due to the factorization property of Section $2$ we have
(for $Q'=Qq^z$):
\[ \calg_K = \lan I_K (Q', q), I_K (Q, q^{-1}) \ran \]
where the series $I_K (Q,q)=p^{\ln Q/\ln q} \sum J_d Q^d $ has coefficients
$J_d\in K^*_{S^1\times G} (X)$. The conjecture we are about to describe 
says that {\em $I_K$ is a common eigen-function of the  
conservation laws of the finite-difference Toda lattice corresponding to
$\Gg'$.} 

\medskip

{\bf 5.2.} The commuting differential operators of quantum Toda lattices 
originate (see \cite{BKo, STS})
 from the center $Z_{U\Gg '}$ of the universal enveloping algebra $U\Gg'$.
Similarly, commuting finite-difference operators of the Toda lattice originate
from the center of the corresponding quantum group $U_q\Gg '$.

 Let $\Gg'$ be the simple complex Lie algebra Langlands-dual to 
$\Gg $ (\this $\Gg'=\Gg$ unless they have the types $B_r$ and $C_r$ which are
dual to one another). The quantum group $U_q\Gg'$ is defined as a (Hopf) 
algebra in terms of generators $K_i^{\pm 1}, X_i^{\pm}, i=1,...,r$ 
and relations (see for instance \cite{RTF}): 
\[ [K_i, K_j]=0, \ K_i X_j^{\pm}=q^{\pm (\a_i, \a_j)} X_j^{\pm} K_i,\ 
[ X_i, X_j ]=\d_{ij} \frac{ K_i-K_i^{-1} }{q-q^{-1}}, \]
\[ \sum_{k=0}^{1-a_{ji}} (-1)^k \binom{1-a_{ji}}{k}_{q^2_i} q_i^{-k(m-k)} 
(X_i^{\pm})^k X_j^{\pm} (X_i^{\pm})^{1-a_{ji}-k} =0 \ \ \text{for}\ i\neq j, \]
where  $a_{ji}=2(\a_i,\a_j)/(\a_i,\a_i)$ is the Cartan matrix of $\Gg'$,
$\a_1,...,\a_r$ --- simple roots of $\Gg'$, $(.,.)$ --- a $W$-invariant
inner product, $q_i=q^{(\a_i,\a_i)/2}$, and $\binom{m}{k}_q$ is the 
$q$-binomial coefficient
\[ \frac{(1-q)(1-q^2)...(1-q^m)}{(1-q)...(1-q^k)\ (1-q)...(1-q^{m-k})}. \]
The Cartan subalgebra $U_q\Gh'=\QQ [K^{\pm 1}]$ is the group algebra of the
root lattice for $\Gg '$. In the Toda theory it is useful to extend the root
lattice (and the quantum group) to the {\em co}-weight lattice. We introduce 
the commuting generators $P_1,...,P_r$ corresponding to
fundamental weights of $\Gg $ so that $K_i=P_1^{(\a_i,\a_1)}...
P_r^{(\a_i,\a_r)}$
and choose new generators $Q_i^{\pm}=X_i^{\pm} 
P_1^{\pm m_{i1}}...P_r^{\pm m_{ir}}$
instead of $X_i^{\pm}$. Here $(m_{ij})$ is any matrix with the property
\[ m_{ij}=m_{ji} \ \text{unless} \ (\a_i,\a_j)<0 \ \text{in which case}\  
m_{ji}-m_{ij}=\pm (\a ,\a )/2 \ , \] 
where $\a $ is a long root. 
(For instance, one can orient the edges of the Dynkin diagram and  
put $m_{ji}=0$ unless there is an edge from $i$ to $j$ in which 
case put $m_{ji}=(\a,\a)/2$.) The relations between the new generators read:
\[ [P_i, P_j]=0, \ P_i Q_j^{\pm}=q^{\pm \d_{ij}} Q_j^{\pm} P_i,\ 
Q_i^{+}Q_j^{-}-q^{m_{ji}-m_{ij}}Q_j^{-}Q_i^{+}=
\d_{ij}\frac{K_i-K_i^{-1}}{q-q^{-1}}, \]
and for $i\neq j$, $x:= q_i^{a_{ji}\pm 2(m_{ji}-m_{ij})/(\a_i,\a_i)}$
\be \sum_{k=0}^{1-a_{ji}}(-1)^k \binom{1-a_{ji}}{k}_{q_i^2} (q_i^2)^{k(k-1)/2}
\ x^k \ (Q_i^{\pm})^k Q_j^{\pm} (Q_i^{\pm})^{1-a_{ji}-k} =0\ .
\label{Serre} \end{equation}
Notice that the Serre relations (\ref{Serre}) allow $1$-dimensional 
representations 
$\x_{-}: U_q\Gn'_{-}$ $\to \CC $ of the subalgebra 
generated by $(Q_1^{-},...,Q_r^{-})$ due to the following
$q$-binomial identity:
\[ (1-x)(1-qx)...(1-q^{m-1}x)=\sum_{k=0}^m (-1)^k \binom{m}{k}_q q^{k(k-1)/2}
x^k\ ,\]
and that the Borel subalgebra $U_q\Gb'_{+}$ generated by $P_i^{\pm 1}$ and 
$Q_j^{+}$ has a representation $\pi_{+}$ onto the algebra of 
finite-difference operators
$Q_i^{+}\mapsto Q_i\times ,\ P_i\mapsto q^{Q_i\p/\p Q_i}, \ i=1,...,r$.

In order to construct commuting finite-difference operators from the center
$Z_{U_q\Gg'}$ of the quantum group we, following B. Kostant, represent the
quantum group as the tensor product $U_q\Gb'_{+}\cdot U_q\Gn'_{-}$ of 
the subspaces, then
notice that the linear map $U_q\Gg'\to U_q\Gb'_{+}\otimes U_q\Gn'_{-}$ 
restricted to the center is a homomorphism of algebras
$Z_{U_q\Gg'} \to U_q\Gb'_{+} \otimes U^{\circ}_q\Gn'_{-}$ (here $\ ^{\circ}$ 
means the anti-isomorphic algebra) and compose the
homomorphism with the representation $\pi_{+}\otimes \x_{-}$.  
   
{\bf 5.3.} 
The center $Z_{U_q\Gg'}$ for generic $q$ is known to be isomorphic to
the polynomial algebra in $r$ generators. The corresponding finite-difference 
operators one obtains by choosing $\x_{-}=0$ are constant coefficient linear 
combinations of translation operators. Considered as Laurent
polynomials of the elementary translations $\hat{P}_i=q^{Q_i\p/\p Q_i}$, 
they become $W$-invariant after the $\rho $-shift 
\[ \hat{P}_i \mapsto \hat{P}_i q^{-\rho_i} = Q^{\rho } \hat{P}_i Q^{-\rho},\] 
where $\sum \rho_i \a_i =\rho $ is
the semi-sum of positive roots of $\Gg'$. The $W$-invariance follows from the
theory of Verma modules for the quantum group $U_q\Gg'$.

One defines commuting operators $\hat{D}_i,\ i=1,...,r$,
of the finite-difference Toda lattice by choosing a generic character on
the role of $\x_{-}$ (\this $\x_{-} (Q_j^{-})\neq 0$ for all $j$) and 
applying the above construction followed by the $\rho$-shift to generators 
$D_1,...,D_r$ of the center. Then the constant coefficients symbols 
$Smb_D:=\hat{D}_i\ \mod (Q)$ 
are $W$-invariant Laurent polynomials in $(\hat{P}_1,...,\hat{P}_r)$. 
Recall now that the generators $P_i$ correspond to fundamental weights of 
$\Gg $. This allows us to consider
the symbols as $W$-invariant functions on the maximal torus of $G$. 
Finally, the form of the finite-difference Toda system we need reads:
\be \hat{D}_i \ I = Smb_{D_i}(\L )\ I, \ \ i=1,...,r, \label{TodaG} 
\end{equation} 
where $\L =(\L_1,...,\L_r)$ are Laurent coordinates on the maximal torus of 
$G$ equal to highest weights of the fundamental representations.   

\medskip

{\bf Conjecture.} {\em The $K^*_{S^1\times G}(X)$-valued formal 
function $I_K(Q, q)$ 
(and respectively the generating function $\calg_K = 
\lan I_K(Qq^z, q), I_K(Q,q^{-1})\ran$  for the flag manifold $X=G/B$ 
satisfies the system (\ref{TodaG}) of finite-difference Toda 
equations corresponding to
the Langlands-dual Lie algebra $\Gg'$.}

\medskip

{\bf 5.4. Remarks.} (1) Operators of the Toda system (\ref{TodaG}) 
depend on the choice of the matrix $(m_{ji}-m_{ij})$ and of the generic 
character 
$\x_{-}$. The latter ambiguity is compensated by rescalings of $Q_1,...,Q_r$.
A canonical choice of the matrix for the classical series $A_r,B_r,
C_r,D_r$ is described in \cite{RTF} (Theorem $12$). At the moment we are not
ready to specify the choice to be made in the above Conjecture. 

(2) The construction of commuting Toda operators is an algebraic version of
the following geometrical description for differential Toda systems: 
interpret the center $Z_{U\Gg'}$ of the universal enveloping algebra as the
algebra of bi-invariant differential operators on the group $G'$ and make
them act on the functions on the maximal torus $T'$ by extending such 
functions to the dense subset $N'_{+} T' N'_{-} \subset G'$ equivariantly
with respect to given generic characters $\x_{\pm}: N'_{\pm}\to \CC^{\times}$
(\this restrict the operators to the sheaf of functions $f$ on $G'$ satisfying 
$f(n_{+} g n^{-1}_{-})=\x_{+} (n_{+}) f(g) \x_{-} (n^{-1}_{-})$). 
The Hamiltonian differential operator $H$ of the quantum Toda system 
is obtained by this construction (followed by the $\rho$-shift) from the 
bi-invariant Laplacian on $G'$. 

(3) By the Harish-Chandra isomorphism, central elements in $U\Gg'$ correspond 
to $ad$-invariant elements in $S(\Gg')$ (or to $W$-invariant polynomials on
the Cartan subalgebra of $\Gg$).  Similarly, central elements in
$U_q\Gg'$ correspond to $W$-invariant functions on the maximal torus of $G$.
This explains why the operator (\ref{Toda}) does not look like a 
finite-difference version of the $2$-nd order differential operator $H$:
the operator $\hat{H}$ corresponds to the trace of unimodular matrices, 
while $H$ corresponds to the trace of square for
traceless matrices. For classical series there are relatively simple algebraic
formulas for generators of $Z_{U_q\Gg'}$ quantizing traces of powers of 
matrices in the vector representation (see \cite{RTF}, Theorem $14$).    
It is not hard to point out explicitly the finite-difference Toda 
operator corresponding in the above Conjecture to the trace in the 
{\em co}-vector representation of a classical group. 

\medskip

{\bf 5.5. Whittaker functions.} In harmonic analysis, one constructs
eigen-functions of the center $Z_{\Gg'}$ as matrix elements of irreducible
representations of $G'$. 
This construction in the case of the Toda system 
(\ref{TodaG}) includes the following ingredients. Let $| v\ran \ \in V$ be a 
{\em Whittaker vector} in an irreducible representation of $U_q\Gg'$,  
\this a common eigen-vector of the elements $Q_i^{-},\ i=1,...r$: 
$Q_i^{-} | v\ran =\x_{-}(Q_i^{-}) |v\ran $. 
Let $\lan a| \ \in V^*$ be a {\em Whittaker covector}: 
$\lan a | Q_i^{+} | x \ran = \x_{+}(Q_i^{+}) \lan a | x \ran $ for any 
$| x \ran \ \in V, \ i=1,...,r$.
Let $Smb_D (\L q^{\rho })$ be the eigenvalue of the central element $D$ in the 
representation $V$. Then the matrix element
\[   \lan a | P_1^{\ln Q_1/\ln q}...\ P_r^{\ln Q_r /\ln q} | v \ran , \]
up to the $\rho$-shift, is an eigenfunction of the finite-difference operator
$\hat{D}$ with the eigenvalue $Smb_D$.

In the case of differential Toda lattices the construction of eigen-functions
as matrix elements (which requires existence of suitable Whittaker vectors,
covectors and integrability of the representation of $\Gg' $ to the maximal
torus $T'\subset G'$ at least) can be realized (see \cite{BKo, STS}) 
in suitable 
analytic versions of the principal series representations on the role of $V$.  
According to \cite{PE} (see the Remark (5) in the Introduction)
the construction can be carried over to the case of quantum groups.
We arrive therefore at the following representation-theoretic interpretation
of the Conjecture 
\[  \begin{array}{c}
\text{generating functions for} \\
\text{representations of}\  S^1\times G \\
\text{in cohomology of line bundles on} \\
\text{spaces of rational curves in}\ G/B \end{array}   \ = \ 
 \begin{array}{c} 
\text{matrix elements in} \\
\text{ $\infty$-dimensional representations} \\
\text{of the Langlands-dual} \\
\text{quantum group} \ U_q\Gg' \end{array}  .\]
We expect that the LHS of this equality has a natural interpretation 
in terms of representation theory for the Kac-Moody loop group $\hat{LG}$ at
the critical level (rather than in terms of generating functions for 
representations of $S^1\times G$). The reason is that the graph spaces 
$GX_d$ (or the quasimap spaces $\mu (GX_d)\subset \Pi_d$) 
can be considered  as degree $d$ approximations 
to the loop space $LX$ with the circle action induced by the rotation of loops.
Moreover, according to \cite{FFKM} these spaces provide adequate geometrical 
background for semi-infinite representation theory of loop groups.   
However we do not have at the moment any direct evidence in favor of this
expectation.

\end{document}